\documentclass[11pt]{article}
\usepackage{epsfig}

\def\be{\begin{equation}}
\def\ee{\end{equation}}
\def\bea{\begin{eqnarray}}
\def\eea{\end{eqnarray}}
\def\bes{\begin{eqnarray*}}
\def\ees{\end{eqnarray*}}

\def\nn{\nonumber}
\def\lb{\label}

\def\T{{\cal T}}
\def\H{{\cal H}}
\def\P{{\cal P}}
\def\G{{\cal G}}
\def\I{{\cal I}}
\def\R{{\bf R}}
\def\C{{\bf C}}
\def\Z{{\bf Z}}

\def\N{{\bf N}}

\def\Q{{\bf Q}}
\def\T{{\bf T}}

\def\Sg{{\Sigma}}

\def\Sg{{\Sigma}}

\def\<{{\langle}}
\def\>{{\rangle}}
\def\T{{\cal T}}

\def\P{{\cal P}}

\def\I{{\cal I}}

\def\Sp{{\rm Sp}}

\def\hb{\vrule height0.18cm width0.14cm $\,$}

\title{On the mean indices of closed characteristics on dynamically convex star-shaped hypersurfaces in $\R^{2n}$}

\author{Wei Wang\thanks{Partially supported by NSFC No. 12025101.
E-mail: wangwei@math.pku.edu.cn  }\\Key Laboratory of Pure and Applied Mathematics\\
School of Mathematical Sciences \\ Peking University, Beijing 100871 \\
PEOPLES REPUBLIC OF CHINA \\ }

\begin{document}
\maketitle

\begin{abstract}
{\it In this paper, we prove that for every dynamically convex compact star-shaped hypersurface
$\Sigma\subset\R^{2n}$,  there exist at least $\lfloor\frac{n+1}{2}\rfloor$ geometrically distinct   closed characteristics  possessing irrational mean indices provided the number of geometrically distinct closed characteristics on $\Sigma$ is finite, this improves Theorem 1.3 in \cite{LoZ} of Y. Long and C. Zhu by finding one more closed characteristic  possessing irrational mean index when $n$ is odd. Moreover, there exist at least   $\lfloor\frac{n+1}{2}\rfloor+1$ geometrically distinct closed characteristics such that the ratio of the mean indices of any two of them is a irrational number provided the number of geometrically distinct  closed characteristics on $\Sigma$ is finite, this improves Theorem 1.2 in \cite{HuO} of X. Hu and Y. Ou when $n$ is odd. In particular, these estimates are sharp for $n=3$. }
\end{abstract}



\renewcommand{\theequation}{\thesection.\arabic{equation}}
\renewcommand{\thefigure}{\thesection.\arabic{figure}}

\setcounter{figure}{0}
\setcounter{equation}{0}
\section{Introduction and main results}

Let $\Sigma\subset \R^{2n}$ be a $C^3$ compact hypersurface strictly star-shaped with respect to the origin.
We denote the set of all such
hypersurfaces by $\H_{st}(2n)$, and denote by $\H_{con}(2n)$ the subset of $\H_{st}(2n)$ consisting of all
strictly convex hypersurfaces. We consider closed characteristics $(\tau, y)$ on $\Sigma$, which are solutions
of the following problem
\be
\left\{\matrix{\dot{y}=JN_\Sigma(y), \cr
               y(\tau)=y(0), \cr }\right. \lb{1.1}\ee
where $J=\left(\matrix{0 &-I_n\cr
        I_n  & 0\cr}\right)$, $I_n$ is the identity matrix in $\R^n$, $\tau>0$, $N_\Sigma(y)$ is the outward
normal vector of $\Sigma$ at $y$ normalized by the condition $N_\Sigma(y)\cdot y=1$. Here $a\cdot b$ denotes
the standard inner product of $a, b\in\R^{2n}$. A closed characteristic $(\tau, y)$ is {\it prime}, if $\tau$
is the minimal period of $y$. Two closed characteristics $(\tau, y)$ and $(\sigma, z)$ are {\it geometrically
distinct}, if $y(\R)\not= z(\R)$. We denote by $\T(\Sigma)$ the set of geometrically distinct
closed characteristics $(\tau, y)$ on $\Sigma\in\H_{st}(2n)$.  A closed characteristic
$(\tau,y)$ is {\it non-degenerate} if $1$ is a Floquet multiplier of $y$ of precisely algebraic multiplicity
$2$; {\it hyperbolic} if $1$ is a double Floquet multiplier of it and all the other Floquet multipliers
are not on ${\bf U}=\{z\in {\bf C}\mid |z|=1\}$, i.e., the unit circle in the complex plane; {\it elliptic}
if all the Floquet multipliers of $y$ are on ${\bf U}$. We call a $\Sigma\in \H_{st}(2n)$ {\it non-degenerate} if
all the closed characteristics on $\Sigma$ are non-degenerate.

There is a long standing conjecture on the
number of closed characteristics on compact convex hypersurfaces in
$\R^{2n}$: \be \,^{\#}\T(\Sg)\ge n, \qquad \forall \; \Sg\in\H_{cov}(2n).
\lb{1.2}\ee
In 1978, P. Rabinowitz in \cite{Rab1} proved
$^\#\T(\Sg)\ge 1$ for any $\Sg\in\H_{st}(2n)$; A. Weinstein in \cite{Wei1}
proved $^\#\T(\Sg)\ge 1$ for any $\Sg\in\H_{con}(2n)$.
When $n\ge 2$,  in 1987-1988, I. Ekeland-L.
Lassoued, I. Ekeland-H. Hofer, and A. Szulkin (cf. \cite{EkL1},
\cite{EkH1}, \cite{Szu1}) proved independently that
$$ \,^{\#}\T(\Sg)\ge 2, \qquad \forall\,\Sg\in\H_{cov}(2n).$$
In \cite{LoZ} of 2002, Y. Long and C. Zhu further proved
\bea\;^{\#}\T(\Sg) \ge\varrho_n(\Sg)\ge\lfloor\frac{n}{2}\rfloor+1, \qquad \forall\, \Sg\in \H_{con}(2n), \lb{1.3}\eea
where $\varrho_n(\Sg)$ is given by Definition 1.1 of \cite{LoZ} and $\lfloor a\rfloor=\max\{k\in\Z\,|\,k\le a\}$ for $a\in\R$.
 In \cite{WHL} of 2007, W. Wang, X. Hu and Y. Long proved $\,^{\#}\T(\Sg)\ge 3$ for
any $\Sg\in\H_{cov}(6)$.  In \cite{Wan3} of 2016, W. Wang proved $\,^{\#}\T(\Sg)\ge \lfloor\frac{n+1}{2}\rfloor+1$ for every $\Sg\in\H_{cov}(2n)$. In \cite{Wan4} of 2016, W. Wang proved $\,^{\#}\T(\Sg)\ge 4$ for every $\Sg\in\H_{cov}(8)$.
In 2024, E. \c{C}ineli, V. Ginzburg, and B. G\"{u}rel in \cite{CGG} proved $\,^{\#}\T(\Sg)\ge n$ for any dynamically convex $\Sg\in\H_{st}(2n)$, which implies conjecture (\ref{1.2}) is true.

For the star-shaped hypersurfaces, in \cite{Gir} of 1984 and \cite{BLMR} of 1985, $\;^{\#}\T(\Sg)\ge n$ for $\Sg\in\H_{st}(2n)$
was proved under some pinching conditions. In \cite{Vit2} of 1989, C. Viterbo proved a generic existence
result for infinitely many closed characteristics on star-shaped hypersurfaces. In \cite{HuL} of 2002, X. Hu and Y. Long proved that $\;^{\#}\T(\Sg)\ge 2$ for any non-degenerate $\Sg\in \H_{st}(2n)$. In \cite{HWZ2} of 2003,
H. Hofer, K. Wysocki and E. Zehnder proved that $\,^{\#}\T(\Sg)=2$ or $\infty$ holds for every non-degenerate
$\Sg\in\H_{st}(4)$ (the same result was proved in \cite{HWZ1} for any $\Sg\in \H_{cov}(4)$) provided that all the stable and unstable manifolds of the  closed characteristics
on $\Sg$ intersect transversally. In 2016, $\;^{\#}\T(\Sg)\ge 2$ was  proved for every $\Sg\in \H_{st}(4)$
by D. Cristofaro-Gardiner and M. Hutchings in \cite{CGH} without any pinching or non-degeneracy conditions,
different proofs of this result can also be found in \cite{GHHM}, \cite{LLo1} and \cite{GiG1}.
In \cite{DLLW}, H. Duan, H. Liu, Y. Long and W. Wang proved $\,^{\#}\T(\Sg)\ge n$
 for every index perfect non-degenerate $\Sg\in \H_{st}(2n)$, and there exist at least $n$ ($n-1$) non-hyperbolic closed characteristics  when $n$ is even (odd). Here $\Sg\in \H_{st}(2n)$ is {\rm index perfect} if it carries only finitely many geometrically distinct prime
closed characteristics, and every prime closed characteristic $(\tau,y)$ on $\Sigma$ possesses positive mean
index and whose Maslov-type index $i(y, m)$ of its $m$-th iteration satisfies $i(y, m)\not= -1$ when $n$ is even,
and $i(y, m)\not\in \{-2,-1,0\}$ when $n$ is odd for all $m\in\N$.
In \cite{GuK}, J. Gutt and J. Kang proved $\,^{\#}\T(\Sg)\ge n$
 for every  non-degenerate $\Sg\in \H_{st}(2n)$ satisfying the condition:
 every
closed characteristic on $\Sg$ possesses Conley-Zehnder index at least $n-1$.
In \cite{GG}, V. Ginzburg and B. G\"{u}rel proved $\,^{\#}\T(\Sg)\ge \lfloor\frac{n+1}{2}\rfloor+1$ for any dynamically convex $\Sg\in\H_{st}(2n)$. In \cite{GGM}, V. Ginzburg, B. G\"{u}rel and L. Macarini studied the multiplicity of closed Reeb orbits on prequantization bundles and obtained
various existence results.  In \cite{AbM}, M. Abreu and L. Macarini studied  closed Reeb orbits for dynamically convex contact forms and obtained various existence and stability results. In 2024, E. \c{C}ineli, V. Ginzburg, and B. G\"{u}rel in \cite{CGG} proved $\,^{\#}\T(\Sg)\ge n$ for any dynamically convex $\Sg\in\H_{st}(2n)$.

For the stability problem,  in
\cite{Eke2} of I. Ekeland in 1986 and \cite{Lon1} of Y. Long in 1998, for
any $\Sg\in\H_{con}(2n)$ the existence of at least one non-hyperbolic
closed characteristic on $\Sg$ was proved provided
$^\#\T(\Sg)<+\infty$. I. Ekeland proved also in \cite{Eke2} the
existence of at least one elliptic closed characteristic on $\Sg$
provided $\Sg\in\H_{con}(2n)$ is $\sqrt{2}$-pinched. In \cite{DDE1} of
1992, Dell'Antonio, D'Onofrio and I. Ekeland proved the existence of at
least one elliptic closed characteristic on $\Sg$ provided
$\Sg\in\H_{con}(2n)$ satisfies $\Sg=-\Sg$. In \cite{LoZ} of 2002, Y. Long
and C. Zhu further proved when $\Sg\in\H_{con}(2n)$ satisfying $^\#\T(\Sg)<+\infty$, there exists
at least one elliptic closed characteristic and there are at least
$\varrho_n(\Sg)-1\ge\lfloor\frac{n}{2}\rfloor$ geometrically distinct closed characteristics on
$\Sg$ possessing irrational mean indices. In \cite{Wan1}, W. Wang proved there are at least
$2$ geometrically distinct closed characteristics possessing irrational mean indices
provided $\Sg\in\H_{con}(6)$ satisfies $^\#\T(\Sg)<+\infty$. In \cite{HuO}, X. Hu and Y. Ou proved the existence of at
least two elliptic closed characteristic on $\Sg$
and
there are at least $\varrho_n(\Sg)\ge\lfloor\frac{n}{2}\rfloor+1$ geometrically distinct closed characteristics  such that the ratio of the mean indices of any two of them is a irrational number
provided
$\Sg\in\H_{con}(2n)$ satisfies $^\#\T(\Sg)<+\infty$.

For the star-shaped case, there are also few results for the stability problem. In \cite{LiL2} of 1999, C. Liu and Y. Long proved that for $\Sg\in\H_{st}(2n)$
either there are infinitely many geometrically distinct closed characteristics, or there exists at least one non-hyperbolic closed
characteristic, provided any closed characteristic possesses its Maslov-type mean index greater than 2
when $n$ is odd, and greater than 1 when $n$ is even. In \cite{LLo2}, H. Liu and Y. Long proved
that $\Sg\in\H_{st}(4)$ and $\,^{\#}\T(\Sg)=2$ imply that both of the closed characteristics must be elliptic
provided that $\Sg$ is symmetric with respect to the origin.

Motivated by the above results, we prove the following results in this paper.

\medskip

{\bf Theorem 1.1.} {\it For every dynamically convex $\Sg\in\H_{st}(2n)$,  there exist at least   $\lfloor\frac{n+1}{2}\rfloor$ geometrically distinct closed characteristics  possessing irrational mean indices provided $^\#\T(\Sg)<+\infty$. }

{\bf Theorem 1.2.} {\it For every dynamically convex $\Sg\in\H_{st}(2n)$,  there exist at least   $\lfloor\frac{n+1}{2}\rfloor+1$ geometrically distinct closed characteristics such that the ratio of the mean indices of any two of them is a irrational number provided $^\#\T(\Sg)<+\infty$. }

{\bf Corollary 1.3.} {\it For every $\Sg\in\H_{con}(2n)$,  there exist at least   $\lfloor\frac{n+1}{2}\rfloor$ geometrically distinct closed characteristics  possessing irrational mean indices provided  $^\#\T(\Sg)<+\infty$. }

{\bf Corollary 1.4.} {\it For every $\Sg\in\H_{con}(2n)$,  there exist at least   $\lfloor\frac{n+1}{2}\rfloor+1$ geometrically distinct closed characteristics  such that the ratio of the mean indices of any two of them is a irrational number provided  $^\#\T(\Sg)<+\infty$. }

{\bf Remark 1.5.} In Theorem 1.3 of \cite{LoZ},  Y. Long and C. Zhu  proved
there are at least $\lfloor\frac{n}{2}\rfloor$ geometrically distinct closed characteristics  possessing irrational mean indices  for $\Sg\in \H_{cov}(2n)$ satisfying
$^{\#}\T(\Sg) <\infty$. Thus Corollary 1.3 obtains one more closed characteristic
possessing irrational mean index than Long-Zhu's result  when
$n$ is odd. In Theorem 1.2 of \cite{HuO},  X. Hu and Y. Ou  proved
there are at least $\lfloor\frac{n}{2}\rfloor+1$ geometrically distinct closed characteristics  such that the ratio of the mean indices of any two of them is a irrational number  for $\Sg\in \H_{cov}(2n)$ satisfying
$^{\#}\T(\Sg) <\infty$. Thus Corollary 1.4 obtains one more closed characteristic
having the same property than Hu-Ou's result  when
$n$ is odd.

{\bf Remark 1.6.} The example of Harmonic Oscillator (cf. Section 1.7 of \cite{Eke3}): Let $0<\alpha_1\le\alpha_2\le\ldots\le\alpha_n$ and
\bea \mathcal{E}_n(r)=\left\{z=(p_1, \ldots,p_n,
q_1,\ldots,q_n)\in\R^{2n}\left |\frac{}{}\right.
\frac{1}{2}\sum_{i=1}^n\alpha_i(p_i^2+q_i^2)=1\right\}.\lb{1.4}\eea
If $\frac{\alpha_i}{\alpha_j}\notin\Q$ whenever $i\neq j$, then
$^{\#}\T(\mathcal{E}_n(r)) = n$ and
$\hat i(x_i)=2\sum_{j=1}^n\frac{\alpha_j}{\alpha_i}$
where $x_i$ is the closed characteristic lying on the $z_i$ plane.

If $\alpha_1=\alpha_2=\cdots =\alpha_n=1$, then $^{\#}\T(\mathcal{E}_n(r)) = \infty$ and
$\hat i(x)=2n$ for any prime closed characteristics, hence there are no closed characteristics
on $\mathcal{E}_n(r)$  possessing irrational mean indices.
Therefore one can not hope to find closed characteristics
 possessing irrational mean indices or a pair of closed characteristics such that the ratio of their mean indices  is a irrational number  when $^{\#}\T(\Sigma) = \infty$
 in general.

 When $^{\#}\T(\Sigma) < \infty$,
 there are $2$ closed characteristics
on $\Sigma$  possessing irrational mean indices for $\Sigma\in\H_{cov}(4)$ (cf. \cite{WHL})
and $\Sigma\in\H_{cov}(6)$ (cf. \cite{Wan1}).
For $n=3$, taking $\alpha_1=2-\sqrt{3},\,\alpha_2=1,\,\alpha_3=\sqrt{3}$,
we have $\frac{\alpha_i}{\alpha_j}\notin\Q$ whenever $i\neq j$
and
\bea &&\hat i(x_1)=\frac{6}{2-\sqrt{3}}\notin\Q,\quad
\hat i(x_2)=6\in\Q,\quad
\hat i(x_3)=\frac{6}{\sqrt{3}}\notin\Q,\quad
\nn\\
&&\frac{\hat i(x_i)}{\hat i(x_j)}=\frac{\alpha_j}{\alpha_i}\notin\Q,\quad 1\le i<j\le 3,\nn
\eea
this example shows that the estimates above are sharp for the case $n=3$.

The proof of Theorem 1.1 is motivated by the methods in
\cite{Gin}, \cite{GG}, \cite{GHHM}, \cite{Hin} and \cite{LoZ}.  In fact, by \cite{LoZ}, there are at least
$\rho_n(\Sigma)\ge\lfloor\frac{n}{2}\rfloor+1$ geometrically distinct
closed characteristics in the common index jump intervals
\bea ^{\#}\T(\Sigma)&\ge& r\nn\\&\ge&
 ^{\#}((2\N-2+n)\cap\cap_{j=1}^r\G_{2m_j-1}(\gamma_{x_j}))
\nn\\&\ge&
 ^{\#}((2\N-2+n)\cap [2N-\kappa_1,\,2N+\kappa_2])\nn\\&\ge&\rho_n(\Sigma)\ge\lfloor\frac{n}{2}\rfloor+1,
 \nn\eea
 where $\kappa_1,\,\kappa_2$ are given by (1.59) and (1.60) of \cite{LoZ}.
 Moreover, $\kappa_2\ge  n-1$ by (5.5) of \cite{LoZ}; $\kappa_1\ge 2$ and  $\kappa_1=2$ if
 and only if there exists a closed characteristic $x$ of the form $P_x=N_1(1, 1)\diamond N_1(1, -1)^{\diamond(n-1)}$ and $i(x, 1)=n$ by (5.15)-(5.18) of \cite{LoZ}. Then there are $\lfloor\frac{n}{2}\rfloor$ elements among these $\lfloor\frac{n}{2}\rfloor+1$ closed characteristics
 possessing irrational mean indices by the Ekeland's variational structure for convex Hamiltonian
 systems in \cite{Eke3}. Now in order to prove Theorem 1.1, we only need to consider the case that $n$ is odd. We use the Lusternik-Schnirelmann theory for the shift operator in equivariant  Floer and symplectic homology developed by V. Ginzburg and B. G\"{u}rel in \cite{GG}.
This theory have almost the same properties as  Ekeland's variational structure.
The advantage for using this structure is that  the  common index jump intervals in this setting contains one more element
comparing with the above one, i.e., it contains $\lfloor\frac{n+1}{2}\rfloor+1$ elements.
In fact, otherwise there must be a  so called symplectic degenerate maximum (SDM), and the
existence of a SDM implies there are infinitely many geometrically distinct closed characteristics (this was proved in \cite{GHHM}, which is motivated by the methods in \cite{Gin} and \cite{Hin} in order to prove the Conley conjecture). Therefore
a similar argument as in \cite{LoZ} moved to this setting shows that there are $\lfloor\frac{n+1}{2}\rfloor$ elements
among these $\lfloor\frac{n+1}{2}\rfloor+1$ closed characteristics possessing irrational mean indices.

The proof of Theorem 1.2 is motivated by the methods in \cite{GG}
and \cite{HuO}. As in the proof of Theorem 1.1, we only need to consider the case that $n$ is odd.
The key ingredient is that we  find that
the  Index recurrence theorem, cf. Theorem 5.2 of \cite{GG} (Theorem 2.4 in this paper)
has a different proof by using the common index jump theorem of Long-Zhu in \cite{LoZ},
and then we have a explicit expression for the variables in Theorem 2.4. Therefore we
can use a similar argument as in \cite{HuO} to obtained the desired result.

In this paper, let $\N$, $\Z$, $\Q$, $\R$ and  $\C$  denote
the sets of natural integers, integers,
rational numbers, real numbers and complex numbers
respectively. Denote by \be \lfloor a\rfloor=\max\{k\in\Z\,|\,k\le
a\}\lb{1.5}\ee
for $a\in\R$.

\setcounter{figure}{0}
\setcounter{equation}{0}
\section{Lusternik-Schnirelmann theory for the shift operator in equivariant Floer and symplectic  homology}

In this section, we review briefly the  Lusternik-Schnirelmann theory for the shift operator in equivariant Floer and symplectic  homology
developed in \cite{GG}, which will be used in the next section in the proof of our main theorems.

Fix a $\Sg\in\H_{st}(2n)$ and assume the following condition holds in the rest of this paper.

(F) {\bf There exist only finitely many geometrically distinct prime closed characteristics\\
$\qquad\qquad \{(\tau_j, x_j)\}_{1\le j\le r}$ (denoted by $\{x_j\}_{1\le j\le r}$ for simplicity in the following) on $\Sigma$. }

Since $\Sg$ is star-shaped, $\Sg$ is differmorphic to the $2n-1$ dimensional sphere
$S^{2n-1}$. As in \S3.2.1 of \cite{GG},
let $\alpha$ be a contact form on $M = S^{2n-1}$ supporting the standard contact structure. Then $(M, d\alpha)$ can be embedded as a hypersurface in $\R^{2n}$ bounding a star-shaped
domain $W$.

Denote by ${\rm SH}^{S^1,+}_\ast(W; \Q)$ the standard positive equivariant symplectic
homology, then we have
\bea {\rm SH}^{S^1,+}_\ast(W; \Q)=\left\{\matrix{&\Q &\quad for \quad \ast= n + 1, n + 3, \ldots, \cr
                  &0& \quad otherwise,\cr}\right. \lb{2.1}
\eea
and the shift operator
\bea D: \Q={\rm SH}^{S^1,+}_\ast(W; \Q)\rightarrow {\rm SH}^{S^1,+}_{\ast-2}(W; \Q)=\Q,\lb{2.2}\eea
which is an isomorphism for $\ast=n+3, n+5,\ldots$.
Moreover there exists a sequence of non-zero homology classes $w_k\in {\rm SH}^{S^1,+}_{n+2k-1}(W; \Q)$ such that $Dw_{k+1}=w_k$ for $k\in\N$. Let $c_k(\alpha)=c_{w_k}(\alpha)$, where $c_{w_k}(\alpha)$ is the spectral invariant for $w_k\in  {\rm SH}^{S^1,+}_{\ast}(W; \Q)$ defined by (3.2) of \cite{GG}.

As in P.546 of \cite{GG},  denote the set of closed Reeb orbits of $\alpha$  by $\P(\alpha)$
and  define the action
\bea A_\alpha(x)=\int_x\alpha,\lb{2.3}\eea
for $x\in\P(\alpha)$.

Now we review the index theory for symplectic paths. We refer the readers to \cite{Lon4}
for  a thorough treatment for this theory. Since we use the structure in \cite{GG} in this paper, it is more convenient to use
the notations in \cite{GG} here.
As is \S4.1 of  \cite{GG}, given any continuous path $\Phi: [0,\,1]\rightarrow\Sp(2m)$
beginning at $\Phi(0)=I_{2m}$, where $\Sp(2m)$ is the symplectic group
defined by $$ \Sp(2m) = \{M\in {\rm GL}(2m,\R)\,|\,M^TJM=J\}, $$
whose topology is induced from that of $\R^{4m^2}$
and $J=\left(\matrix{0 &-I_m\cr
        I_m  & 0\cr}\right)$, the upper
and lower Conley-Zehnder indices $\mu_+(\Phi)$, $\mu_-(\Phi)$  and the mean index $\hat\mu(\Phi)$ can be defined, cf. \cite{Lon4} for details.
In particular, we have
\bea \hat\mu(\Phi^k)=k\hat\mu(\Phi),\quad k\in\N.\nn\eea

Let $(M,\,\alpha)$ be the contact manifold as above and denote by $\varphi_t$ its Reeb flow.
For a closed Reeb orbit (closed characteristic) $x$, its linearized Poincar\'e
map $\Phi=d\varphi_t|_x\in\tilde\Sp(2m)$ (where $\tilde{\Sp}(2m)$ is the universal covering of $\Sp(2m)$) with $m=n-1$ is defined in the standard way, cf. P.561 of \cite{GG}. Then $\mu_\pm(x), \hat\mu(x)$, etc. are defined by $\mu_\pm(\Phi), \hat\mu(\Phi)$, etc. In particular, we have
\bea \hat\mu(x^k)=k\hat\mu(x),\quad k\in\N,\lb{2.4}\eea
where $x^k$ is the k-th iteration of $x$ defined by $x^k(t)=x(kt)$.

Comparing to Theorems 1.5 and 1.6
of \cite{LoZ} and Proposition 3.5 of \cite{Wan1}, we have the following  Floer and symplectic version of these results.

{\bf Proposition 2.1.} (cf. Corollary 3.9 of \cite{GG}){\it Under the assumption (F), we have
\bea c_1(\alpha)<c_2(\alpha)<c_3(\alpha)<\cdots.\nn\eea
There exists an injection
$$\phi:\N\rightarrow\P(\alpha), \quad k\mapsto y_k$$
called a carrier map such that $c_k(\alpha)=A_\alpha(y_k)$ and
\be A_\alpha(y_1)<A_\alpha(y_2)<A_\alpha(y_3)<\cdots,\lb{2.5}\ee
and  ${\rm SH}^{S^1}_{n+2k-1}(y_k; \,\Q)\neq 0$,
where ${\rm SH}^{S^1}_\ast(y_k; \,\Q)$ is the equivariant local symplectic homology of $y_k$
as in Definition 3.5 of \cite{GG}.  Moreover,
\bea \mu_-(y_k)\le n+2k-1\le \mu_+(y_k),\lb{2.6}\\
|\hat\mu(y_k)-(n+2k-1)|\le n-1.\lb{2.7}\eea
}

The  following definition of {\it reoccurring}  in \S6.1.1 of \cite{GG} is the  Floer and symplectic
version of the definition of {\it infinite variationally visible}  in Definition 1.4 of \cite{LoZ}.

{\bf Definition 2.2.} {\it A prime closed characteristic $x$ is reoccurring if its iterations occur infinitely many times in the image of a carrier injection $\phi$ from Proposition 2.1.
}

For $x\in\P(\alpha)$. denote by
\be \hat c(x)=\frac{A_\alpha(x)}{\hat\mu(x)}.\lb{2.8}\ee
Then by (\ref{2.3}) and (\ref{2.4}),
\be \hat c(x^k)=\frac{A_\alpha(x^k)}{\hat\mu(x^k)}
=\frac{\int_{x^k}\alpha}{\hat\mu(x^k)}=\frac{k\int_{x}\alpha}{k\hat\mu(x)}
=\frac{\int_{x}\alpha}{\hat\mu(x)}=\hat c(x).\lb{2.9}\ee

By Theorem 6.4 of \cite{GG}, we have the following resonance relations, this is
the  Floer and symplectic version of Theorem 5.3.15 of \cite{Eke3}.

{\bf Proposition 2.3.} {\it For any two reoccurring closed characteristics $x$ and $y$
we have $\hat c(x)=\hat c(y)$, i.e.,
\be  \frac{A_\alpha(x)}{\hat\mu(x)}=\frac{A_\alpha(y)}{\hat\mu(y)} .\lb{2.10}\ee
}

We have the following Index recurrence theorem, cf. Theorem 5.2 of \cite{GG}. This theorem
 plays the similar role as the common index jump theorem of Long-Zhu, cf. Theorem 4.3
in \cite{LoZ} and its enhanced version, cf.  Theorem 3.5 in \cite{DLW}.

 {\bf Theorem 2.4.} {\it Let $\Phi_1,\ldots, \Phi_r$ be a finite collection of
elements in $\tilde{\Sp}(2m)$  satisfying $\hat\mu(\Phi_i)>0$ for $1\le i\le r$. Then for any
$\eta>0$ and any $l_0\in\N$,
there exists an integer sequence $d_j\rightarrow\infty$ and $r$ integer sequences
$k_{ij}\rightarrow \infty$ for $1\le i\le r$ such that for $1\le i\le r$, $j\in\N$ and $1\le l\le l_0$,
we have
\bea &&|\hat\mu(\Phi_i^{k_{ij}})-d_j|<\eta\lb{2.11}\\
&&\mu_\pm(\Phi_i^{k_{ij}+l})=d_j+\mu_\pm(\Phi_i^l),\lb{2.12}\\
&&\mu_\pm(\Phi_i^{k_{ij}-l})=d_j+\mu_\pm(\Phi_i^{-l})+b_+(\Phi_i^l)-b_-(\Phi_i^l),\lb{2.13}\eea
where $b_\ast(\Phi)$ is the signature multiplicities of $\Phi$ defined in Definition 4.5 of \cite{GG}.
Moreover, for any $N\in\N$, we can make all $d_j$ and $k_{ij}$ divisible by $N$.
}

As in Definitions 4.6 and 4.9 of \cite{GG}, we make the following definition.

{\bf Definition 2.5.} {\it A path $\Phi\in\tilde\Sp(2m)$ is dynamically convex if $\mu_-(\Phi)\ge m+2$. The Reeb flow on a $(2n-1)$-dimensional contact manifold is said to be dynamically convex if every closed Reeb orbit $x$ is dynamically convex, i.e., $\mu_-(x)\ge  m+2=n+1$.
}

In particular, we have

{\bf Theorem 2.6.} (cf. Theorem 4.10 of \cite{GG}) {\it The Reeb flow on a strictly convex hypersurface in $\R^{2n}$
is dynamically convex.
}

Under the assumption of dynamically convexity, we have the following estimates
by Corollary 5.4 of \cite{GG}.

{\bf Proposition 2.7.} {\it Suppose the paths $\Phi_1,\ldots,\Phi_r$ in Theorem 2.4
are dynamically convex. Then for $1\le l\le l_0$,
\bea &&\mu_-(\Phi_i^{k_{ij}+l})\ge d_j+2l+m,\lb{2.14}\\
&&\mu_+(\Phi_i^{k_{ij}-l})\le d_j-2l.\lb{2.15}\eea
In particular, $\mu_-(\Phi_i^{k_{ij}+l})\ge d_j+2+m$ and
$\mu_+(\Phi_i^{k_{ij}-l})\le d_j-2$ for all $l\in\N$.
}

\setcounter{figure}{0}
\setcounter{equation}{0}
\section{Proof of the main theorems}

In his section, we  give the proof of the main theorems. The proof of Theorem 1.1 is motivated by the methods
of Theorem 1.3 in \cite{LoZ}, the  difference here is  that we work in the Floer and symplectic framework,
and then the common index jump intervals has one more element than that of \cite{LoZ}
when $n$ is odd.
The new property that the common index jump intervals here has one more element
is proved by contradiction: suppose the contrary, then there must be a  so called symplectic degenerate maximum (SDM) in the sense of \cite{GHHM}, and the
existence of a SDM implies there must be infinitely many geometrically distinct closed characteristics,
contradicts to the assumption (F).

{\bf Proof of Theorem 1.1.}  Let $\Phi_i\in\tilde{\Sp}(2m)$ with $m=n-1$ be the linearized
Poincar\'{e} map along $x_i$ for $1\le i\le r$. By the assumption of dynamically convexity,
we have
\bea \mu_-(x_i^{k_{ij}+l})\ge d_j+n+1,\quad \mu_+(x_i^{k_{ij}-l})\le d_j-2,
\qquad \forall l, j\in\N,\;1\le i\le r,\lb{3.1}\eea
by Proposition 2.7.

As in P.569  \cite{GG}, Let
\be \psi: \I=\{n+1, n+3,n+5,\ldots\}\rightarrow\P(\alpha),\quad d\mapsto y_d\lb{3.2}\ee
be the composition of the map $\phi$ in Proposition 2.1 with the bijection $\I\rightarrow \N: d\mapsto(d+1-n)/2$, then we have
\bea A_\alpha(y_{n+1})<A_\alpha(y_{n+3})<A_\alpha(y_{n+5})<\cdots,\lb{3.3}\eea
and  ${\rm SH}^{S^1}_d(y_d; \,\Q)\neq 0$   and
\bea \mu_-(y_d)\le d\le \mu_+(y_d),\lb{3.4}\\
|\hat\mu(y_d)-d|\le n-1.\lb{3.5}\eea

Denote by $L=[d_j-1,\,d_j+n]\cap\I$, where $d_j$ is the integer found in Theorem 2.4 for some $j\in\N$ large enough.
By Proposition 2.1 and (\ref{3.1}), for any $d\in L$, $y_d$ must has the form $x_{i(d)}^{k_{i(d)j}}$ for some $i(d)\in\{1,2,\ldots,r\}$.
In fact, $\mu_-(\Phi_i^{k_{ij}+l})\ge d_j+n+1$ implies
${\rm SH}^{S^1}_{d}(x_i^{k_{ij}+l}; \,\Q)=0$ for $1\le i\le r$ and $l\in\N$;
$\mu_+(\Phi_i^{k_{ij}-l})\le d_j-2$ implies
${\rm SH}^{S^1}_{d}(x_i^{k_{ij}-l}; \,\Q)=0$ for $1\le i\le r$ and $l\in\N$.
Moreover, for any $d_1\neq d_2\in L$, we have $i(d_1)\neq i(d_2)$, this follows from
the fact that $A_\alpha(y_{d_1})\neq A_\alpha(y_{d_2})$ together with
$y_{d_1}=x_{i(d_1)}^{k_{i(d_1)j}}$ and $y_{d_2}=x_{i(d_2)}^{k_{i(d_2)j}}$.

Now we consider the case that $n$ being odd.
By Proposition 2.7, $y_{d_j-2}$ can only has the form $x_i^{k_{ij}}$
or $x_i^{k_{ij}-1}$ for some $i\in\{1,2,\ldots,r\}$.  In fact, $\mu_-(\Phi_i^{k_{ij}+l})\ge d_j+2+m$ implies
${\rm SH}^{S^1}_{d_j-2}(x_i^{k_{ij}+l}; \,\Q)=0$ for $1\le i\le r$ and $l\in\N$;
$\mu_+(\Phi_i^{k_{ij}-l})\le d_j-2l$ implies
${\rm SH}^{S^1}_{d_j-2}(x_i^{k_{ij}-l}; \,\Q)=0$ for $1\le i\le r$ and $l\ge 2$.

 If $y_{d_j-2}=x_i^{k_{ij}-1}$, then by P.570 of \cite{GG},
$x_i$ must be a symplectic degenerate maximum (SDM), then there are infinitely many
geometrically distinct closed Reeb orbits (this was proved in \cite{GHHM}, which is motivated by the methods in \cite{Gin} and \cite{Hin} for the proof of the Conley conjecture), contradict to the assumption (F), therefore $y_{d_j-2}=x_i^{k_{ij}}$ for some $i\in\{1, 2,\ldots,r\}$.

Denote by $\tilde L=[d_j-2,\,d_j+n]\cap\I$ for $n$ being odd and $\tilde L=[d_j-1,\,d_j+n]\cap\I$ for $n$ being even, then by the above argument, for any $d\in \tilde L$, $y_d$ must has the form $x_{i(d)}^{k_{i(d)j}}$ for some $i(d)\in\{1,2,\ldots,r\}$; moreover, for any $d_1\neq d_2\in \tilde L$, we have $i(d_1)\neq i(d_2)$.
Note that we can require that $d_j$ being even by Theorem 2.4, therefore we have
the following:

{\bf Claim 1.} {\it There are at least $^\#(\tilde L)=\lfloor \frac{n+1}{2}\rfloor+1\equiv p$ geometrically distinct
closed Reeb orbits.}

Note that the following claim implies Theorem 1.1.

{\bf Claim 2.} {\it Among the $p$ closed Reeb orbits found in Claim 1, there is at  most one
possessing rational mean index.}

We prove the claim by contradiction. Assume
$d_1\neq d_2\in \tilde L$, such that $\hat\mu(x_{i(d_1)}), \hat\mu(x_{i(d_2)})\in\Q$.
Then we can require that the number $N$ in Theorem 2.4 satisfies
$N\hat\mu(x_{i(d_1)}),  N\hat\mu(x_{i(d_2)})\in\Z$.
Since $N|k_{ij}$, we have
\bea\hat\mu(x_{i(d_1)}^{k_{i(d_1)j}})=k_{i(d_1)j}\hat\mu(x_{i(d_1)})
\in\Z,\quad \hat\mu(x_{i(d_2)}^{k_{i(d_2)j}})=k_{i(d_2)j}\hat\mu(x_{i(d_2)})
\in\Z,\lb{3.6}\eea
Therefore by (\ref{2.11}), we have
\bea\hat\mu(x_{i(d_1)}^{k_{i(d_1)j}})=d_j=\hat\mu(x_{i(d_2)}^{k_{i(d_2)j}}).\lb{3.7}\eea

Since the number of prime closed Reeb orbits is finite, we can choose
$x_d$  to be reoccurring for  all $d\in\tilde L$, cf. Theorem 6.1 of \cite{GG}.
Then by (\ref{2.9}) and (\ref{2.10}), we have
\be  \frac{A_\alpha(x_{i(d_1)}^{k_{i(d_1)j}})}{\hat\mu(x_{i(d_1)}^{k_{i(d_1)j}})}=
\frac{A_\alpha(x_{i(d_1)})}{\hat\mu(x_{i(d_1)})}=
\frac{A_\alpha(x_{i(d_2)})}{\hat\mu(x_{i(d_2)})}=\frac{A_\alpha(x_{i(d_2)}^{k_{i(d_2)j}})}{\hat\mu(x_{i(d_2)}^{k_{i(d_2)j}})} .\lb{3.8}\ee
Thus we have
\be A_\alpha(x_{i(d_1)}^{k_{i(d_1)j}})=A_\alpha(x_{i(d_2)}^{k_{i(d_2)j}}).\lb{3.9}\ee
This contradicts to (\ref{3.3}). The proof of Theorem 1.1 is complete.\hfill\hb

{\bf Proof of Theorem 1.2.}  We continue to use the symbols in \cite{GG}, these symbols
are somewhat different than the symbols used in \cite{LoZ}, and we will mention the difference whenever used below.

 Denote by $\exp(\pm2\pi\sqrt{-1}\lambda_{iq})$ the elliptic eigenvalues of $\Phi_i$ with
 irrational $\lambda_{iq}$ for  $1\le i\le r$, $1\le q\le \iota_i$, and $\Delta_i=\hat\mu(\Phi_i)$ for $1\le i\le r$. Then by P.568 of \cite{GG}, in oder to prove Theorem 2.4,
 it is sufficient to find an infinite sequence of tuples $(k_1, k_2,\ldots,k_r)$ such that
 \bea  ||k_i\lambda_{iq}||&<&\epsilon,\quad 1\le i\le r,\;1\le q\le\iota_i,\lb{3.10}\\
 |k_1\Delta_1-k_i\Delta_i|&<&\frac{1}{8}, \quad 2\le i\le r,\lb{3.11}
 \eea
for any given $\epsilon>0$, where $||\cdot||$ stands for the distance to the nearest integer, cf. P.566 of \cite{GG}.
The original proof of (\ref{3.10}) and (\ref{3.11}) in \cite{GG} uses Minkowski's theorem to obtain $(k_1, k_2,\ldots,k_r)$ abstractly. Now we use the method in \cite{LoZ} to construct $(k_1, k_2,\ldots,k_r)$ concretely.

Firstly by (4.10), (4.40) and (4.41) in \cite{LoZ} , let
\bea k_i=\left(\left\lfloor\frac{T}{M\Delta_i}\right\rfloor+\chi_i\right)M,\quad 1\le i\le r,\lb{3.12}\eea
where $M\in\N$ is any fixed integer divisible by $N$ in Theorem 2.4,
$\chi_i=0$ or $1$ and $T\in\N$ will be determined later.
Note that in \cite{LoZ}, they use the symbol $[a]$ stands for $\lfloor a\rfloor$
and $\{a\}$ stands for $a-\lfloor a\rfloor$ for $a\in\R$. In the following, we also use the symbol
$\{a\}=a-\lfloor a\rfloor$.

As in (4.21) of \cite{LoZ},  let $h=r+\sum_{1\le i\le
r}\iota_i$ and \bea  v=\left(\frac{1}{M\Delta_1},\dots,
\frac{1}{M\Delta_r}, \frac{2\lambda_{11}}{\Delta_1}, \frac{2\lambda_{12}}{\Delta_1},
\dots\frac{2\lambda_{1\iota_1}}{\Delta_1}, \frac{2\lambda_{21}}{\Delta_2},\dots, \frac{2\lambda_{r\iota_r}}{\Delta_r}\right)\in\R^h. \lb{3.13}\eea

{\bf Claim 3.} { \it In order to make (\ref{3.10}) and (\ref{3.11}) hold, it is sufficient to find a vertex
\bea\chi=(\chi_1,\ldots,\chi_r,\chi_{11},\chi_{12},\ldots,\chi_{1\iota_1}, \chi_{21},\ldots,\chi_{r\iota_r})\lb{3.14}\eea
of the cube $[0,\,1]^h$ and infinitely many
integers $T\in\N$ such that \bea|\{Tv\}-\chi|<\epsilon\lb{3.15}\eea
for any given $\epsilon>0$ small enough.}

In fact, (\ref{3.10}) follows directly from (4.5) of \cite{LoZ}. For (\ref{3.11}), we have
\bea &&|k_1\Delta_1-k_i\Delta_i|\nn\\=&&\left|\left(\left\lfloor\frac{T}{M\Delta_1}\right\rfloor+\chi_1\right)M\Delta_1-\left(\left\lfloor\frac{T}{M\Delta_i}\right\rfloor+\chi_i\right)M\Delta_i\right|\nn\\=&&\left|\left(\frac{T}{M\Delta_1}-\left\{\frac{T}{M\Delta_1}\right\}+\chi_1\right)M\Delta_1-\left(\frac{T}{M\Delta_i}-\left\{\frac{T}{M\Delta_i}\right\}+\chi_i\right)M\Delta_i\right|\nn\\=&&\left|\left(-\left\{\frac{T}{M\Delta_1}\right\}+\chi_1\right)M\Delta_1-\left(-\left\{\frac{T}{M\Delta_i}\right\}+\chi_i\right)M\Delta_i\right|\nn\\
\le&& 2M\max_{1\le i\le r}|\Delta_i|\epsilon,\lb{3.16}\eea
by (\ref{3.15}).  Therefore (\ref{3.11}) follows by taking
$\epsilon<\frac{1}{16M\max_{1\le i\le r}|\Delta_i|}$.

The vertex $\chi$ in (\ref{3.15}) can be obtained in the following way
(cf. Theorem 4.2 of \cite{LoZ}).  Let $H$ be the closure of $\{\{mv\}|m\in\N\}$ in $\bf {T}^h=(\R/\Z)^h$ and
$V=T_0\pi^{-1}H$ be the tangent space of $\pi^{-1}H$ at the origin in $\R^h$,
where $\pi: \R^h\rightarrow  \bf{T}^h$ is the projection map.
Define
\bea A(v)=V\setminus\cup_{v_k\in\R\setminus\Q}\{x=(x_1,,\ldots,x_h)\in V | x_k=0\}.
\lb{3.17}\eea
Define $\psi(x)=0$  when $x\ge 0$ and $\psi(x)=1$ when $x<0$.
Then for any $a=(a_1,\ldots, a_h)\in A(V)$, the vector
\bea \chi=(\psi(a_1),\ldots,\psi(a_h))\lb{3.18}\eea
makes (\ref{3.15}) holds for infinitely many $T\in\N$.
Moreover, the set $A(v)$ possesses the following properties:

(i) When $v\in\R^h\setminus\Q^h$,  then $\dim V\ge 1$,
$0\notin A(v)\subset V$, $A(v)=-A(v)$ and $A(v)$ is open in $V$.

(ii) When $\dim V = 1$, then  $A(v) = V \setminus\{0\}$.

(iii) When $\dim V \ge 2$,  $A(v)$  is obtained from $V$ by deleting all the
coordinate hyperplanes  with    dimension
strictly smaller than   $\dim V$    from    $V$.

Now Let $T_j\in\N$   satisfying (\ref{3.15})  for some $j\in\N$ large enough,
\bea k_{ij}=\left(\left\lfloor\frac{T_j}{M\Delta_i}\right\rfloor+\chi_i\right)M
,\quad 1\le i\le r, \nn\eea
and $d_j=[k_{ij}\Delta_i]$, where $[\,\cdot\,]$ denotes the nearest integer, cf. P. 567-578 of
\cite{GG}.
Let $\tilde L$ be the set in Claim 1 above,  then $^\#(\tilde L)=\lfloor \frac{n+1}{2}\rfloor+1\equiv p$.

{\bf Claim 4.} {\it Any two  distinct closed Reeb orbits $x_{i(d_1)}$ and $x_{i(d_2)}$ among the $p$ closed Reeb orbits found
 in Claim 1 satisfies $\frac{\hat\mu(x_{i(d_1)})}{\hat\mu(x_{i(d_2)})}\notin\Q$.}

Note that by the dynamically convexity assumption,
we have $\hat\mu(x_i)\ge 2$ for $1\le i\le r$ by Lemma 4.8 of \cite{GG}. We prove  Claim 4 by contradiction. Assume
$d_1\neq d_2\in \tilde L$ such that
\bea\frac{\hat\mu(x_{i(d_1)})}{\hat\mu(x_{i(d_2)})}=\frac{p}{q}\in\Q,\lb{3.19}\eea
where $p, q\in\N$ and $(p, q)=1$. If $\hat\mu(x_{i(d_1)})\in\Q$,
then $\hat\mu(x_{i(d_2)})\in\Q$ also by (\ref{3.19}), this contradicts to
Claim 1 above. Therefore $\hat\mu(x_{i(d_1)}),\,\hat\mu(x_{i(d_2)})\notin\Q$.

By Theorem 1.1, there are at least $\lfloor\frac{n+1}{2}\rfloor$ closed Reeb orbits possessing
irrationally mean indices, therefore $v\in\R^h\setminus\Q^h$ and then $\dim A(v)\ge 1$.
Thus we can choose $a\in A(v)$  small enough and  $T_j\in \N$ such that
$\{T_jv\}-\chi(a)\in V$ and (\ref{3.15}) holds.
By (\ref{3.13}) and $\hat\mu(x_i)=\Delta_i$,  we have
\bea \frac{v_{i(d_1)}}{v_{i(d_2)}}=\frac{\frac{1}{M\Delta_{i(d_1)}}}{\frac{1}{M\Delta_{i(d_2)}}}=\frac{\Delta_{i(d_2)}}{\Delta_{i(d_1)}}=\frac{\hat\mu(x_{i(d_2)})}{\hat\mu(x_{i(d_1)})}
=\frac{q}{p}>0,\lb{3.20}\eea
therefore by the definition of $H=\overline{\{\{mv\}|m\in\N\}}$ in $\bf {T}^h=(\R/\Z)^h$ and
$V=T_0\pi^{-1}H$, we must have
$\frac{x_{i(d_1)}}{x_{i(d_2)}}=\frac{q}{p}>0$ for any $x\in V$ with $x_{i(d_2)}\neq 0$.
In particular, if $a\in A(v)$, then $a_{i(d_1)}>0$, $a_{i(d_2)}>0$  or $a_{i(d_1)}<0$, $a_{i(d_2)}<0$,
hence $\chi_{i(d_1)}=\chi_{i(d_2)}=0$ or $\chi_{i(d_1)}=\chi_{i(d_2)}=1$, i.e.,
$\chi_{i(d_1)}=\chi_{i(d_2)}$ for any $a\in A(v)$.

Since $\{T_jv\}-\chi(a)\in V$, we have
\bea \frac{\{T_jv_{i(d_1)}\}-\chi_{i(d_1)}(a)}{\{T_jv_{i(d_2)}\}-\chi_{i(d_2)}(a)}
=\frac{\{\frac{T_j}{M\Delta_{i(d_1)}}\}-\chi_{i(d_1)}(a)}{\{\frac{T_j}{M\Delta_{i(d_2)}}\}-\chi_{i(d_2)}(a)}=\frac{q}{p}.\lb{3.21}
\eea
From (\ref{3.20}) and (\ref{3.21}), we have
\bea 0&=&\left(\left\{\frac{T_j}{M\Delta_{i(d_1)}}\right\}-\chi_{i(d_1)}(a)\right)M\Delta_{i(d_1)}
-\left(\left\{\frac{T_j}{M\Delta_{i(d_2)}}\right\}-\chi_{i(d_2)}(a)\right)M\Delta_{i(d_2)}
\nn\\&=&\left(\frac{T_j}{M\Delta_{i(d_1)}}-\left\lfloor\frac{T_j}{M\Delta_{i(d_1)}}\right\rfloor-\chi_{i(d_1)}(a)\right)M\Delta_{i(d_1)}
\nn\\&&-\left(\frac{T_j}{M\Delta_{i(d_2)}}-\left\lfloor\frac{T_j}{M\Delta_{i(d_2)}}\right\rfloor-\chi_{i(d_2)}(a)\right)M\Delta_{i(d_2)}
\nn\\&=&-\left(\left\lfloor\frac{T_j}{M\Delta_{i(d_1)}}\right\rfloor+\chi_{i(d_1)}(a)\right)M\Delta_{i(d_1)}
+\left(\left\lfloor\frac{T_j}{M\Delta_{i(d_2)}}\right\rfloor+\chi_{i(d_2)}(a)\right)M\Delta_{i(d_2)}
\nn\\&=&-k_{i(d_1)j}\Delta_{i(d_1)}+k_{i(d_2)j}\Delta_{i(d_2)}\nn\\
&=&-k_{i(d_1)j}\hat\mu(x_{i(d_1)})+k_{i(d_2)j}\hat\mu(x_{i(d_2)})\nn\\
&=&-\hat\mu(x_{i(d_1)}^{k_{i(d_1)j}})+\hat\mu(x_{i(d_2)}^{k_{i(d_2)j}}),\lb{3.22}
\eea
i.e., we have
\bea \hat\mu(x_{i(d_1)}^{k_{i(d_1)j}})=\hat\mu(x_{i(d_2)}^{k_{i(d_2)j}}).\lb{3.23}\eea
Then as in the proof of Theorem 1.1, we have
\bea A_\alpha(x_{i(d_1)}^{k_{i(d_1)j}})=A_\alpha(x_{i(d_2)}^{k_{i(d_2)j}}).\lb{3.24}
\eea
This contradicts to (\ref{3.3}). The proof of Theorem 1.2 is complete.\hfill\hb

{\bf Proof of Corollaries 1.3 and 1.4.} By Theorem 2.6, the Reeb flow on a strictly convex hypersurface in $\R^{2n}$ is dynamically convex, hence Corollaries 1.3 and 1.4 follows from Theorems 1.1 and 1.2.\hfill\hb

\noindent {\bf Acknowledgements.} I would like to sincerely thank my
Ph. D. thesis advisor, Professor Yiming Long, for introducing me to Hamiltonian
dynamics and for his valuable help and encouragement during my research.

\bibliographystyle{abbrv}

\end{document}